\documentclass[11pt]{article}
\usepackage{graphicx}
\usepackage{amsmath,amssymb,amsthm,amsfonts
}
\usepackage{amssymb}
\newtheorem{thm}{Theorem}[section]

\newtheorem{lem}{Lemma}[section]

\newtheorem{rem}{Remark}[section]
\theoremstyle{definition}

\numberwithin{equation}{section}

\DeclareMathSymbol{\C}{\mathalpha}{AMSb}{"43}

\newcommand{\eps}{\varepsilon}

\newcommand{\R}{{\mathbb{R}}}

\newcommand{\inte}{\int_{\mathbb{R}^3}}

\newcommand{\bsub}{\begin{subequations}}
\newcommand{\esub}{\end{subequations}$\!$}

\begin{document}

\title{ Cylindrically Symmetric Ground State Solutions for Curl-Curl Equations with Critical  Exponent \thanks{Email: zengxy09@126.com}}

\author{\ Xiaoyu Zeng\\
\small \it Department of Mathematics, School of Sciences, Wuhan University of Technology,\\
\small \it Wuhan 430070, P. R. China.
  }

\date{\today}

\smallbreak \maketitle

\begin{abstract}
We study the following   nonlinear critical curl-curl equation \begin{equation}\label{eq0.1}\nabla\times \nabla\times U +V(x)U=|U|^{p-2}U+ |U|^4U,\quad  x\in
\mathbb{R}^3,\end{equation}
where $V(x)=V(r, x_3)$ with $r=\sqrt{x_1^2+x_2^2}$ is  1-periodic in $x_3$ direction and belongs to $L^\infty(\R^3)$. When   $0\not\in \sigma(-\Delta+\frac{1}{r^2}+V)$ and $p\in(4,6)$, we prove the existence of nontrivial solution for (\ref{eq0.1}), which is indeed a  ground state solution in a suitable cylindrically symmetric space. Especially, if $ \sigma(-\Delta+\frac{1}{r^2}+V)>0$, a ground state solution is obtained for any $p\in(2,6)$.
\end{abstract}

\vskip 0.2truein


\vskip 0.2truein

\section{Introduction}

In this paper, we consider the following nonlinear critical curl-curl equation
\begin{equation}\label{d1.1}
\nabla\times \nabla\times U +V(x)U=|U|^{p-2}U+ |U|^4U, \quad  x\in
\mathbb{R}^3,
\end{equation}
where $V(x)\in L^{\infty}(\R^3)$ is the external potential,  $p\in(2,6)$ is of subcritical growth. Equation (\ref{d1.1}) is related to the standing waves of nonlinear Maxwell's equation for an inhomogeneous medium. For more physical backgrounds one can refer to \cite{BD,BM,BF}, etc.

 Solution $U\in H({\rm curl}; \R^3)$ of (\ref{d1.1}) is in general a critical point of the energy functional
\begin{equation*}
 I(U)=\frac{1}{2}\inte|\nabla \times U|^2+V(x)|U|^2dx-\frac{1}{p}\inte|U|^pdx-\frac{1}{6}\inte|U|^6dx.
\end{equation*}
The main difficulty for dealing with equation (\ref{d1.1}) lies in the presence of the curl-curl operator, which causes the energy functional $I$   strong indefinite in some sense. To overcome this difficulty,  much research is devoted to the following equation
\begin{equation}\label{eq1.2}
\nabla\times \nabla\times U +V(x)U=g(x,U), \quad  x\in
\mathbb{R}^3,
\end{equation}
where $g(x,U)$ is a subcritical nonlinearity.
For the case of $V(x)\equiv0$, Benci, Fortunato \cite{BF} and Azzollini et al. \cite{ABD} studied equation (\ref{eq1.2}) in a cylindrically symmetric vector space,  which  is instructed with   divergence free elements. In a different symmetric space D'April and Siciliano also obtained in \cite{DS} nontrivial solutions for (\ref{eq1.2}). When $V(x)\not\equiv0$ is cylindrically symmetric, solutions with some kind of symmetries  were also  studied in \cite{BD,HR}. Recently, Bartsch and Mederski \cite{BM} developed the Nehari-Pankov method \cite{SW} and considered the ground state and bound state of (\ref{eq1.2}) in  bounded domains with the boundary condition $\nu\times U=0$. Based on this work, Mederski \cite{M} further studied (\ref{eq1.2}) where, e.g., $g\sim |U|^{q-1}U$ if $|U|\ll1$ and $g\sim |U|^{q-1}U$ if $|U|\gg1$ for $1<p<5<q$ and  $V(x)\leq0$ is periodic and $V\in L^\frac{p-1}{p-1}(\R^3)\cap L^\frac{q-1}{q-1}(\R^3)$. One of the main contributions of \cite{BM,M} is that they can treat equation (\ref{eq1.2}) without any assumption on the symmetry of $U$. For the case that   $g$ is asymptotically linear growth at infinity, we would like to mention some existence results in  \cite{QT} by Qin and Tang, and  some related topics in  \cite{S,SZ1,SZ2,SZ3} by Stuart and Zhou.

Very recently,   Mederski considered   the  Brezis-Nirenberg problem for equation (\ref{d1.1}) (i.e., removing the term $|U|^{p-2}U$) in a bounded domain, and obtained some  existence and no-existence results of cylindrically  symmetric  ground state solutions  for (\ref{d1.1}) in \cite{M2}.  In this paper,  we focus on  problem (\ref{d1.1}) in  a suitable cylindrically symmetric space. We always  assume that  $V(x)$ satisfies
\begin{equation*}
V(x)=V(r,x_3)\in L^\infty(\R^3), \text{ where }r=\sqrt{x_1^2+x_2^2}.
\end{equation*}
One can search for solutions of the form \cite{ABD,BD}
\begin{equation}\label{eq1.04}
U(x)=\frac{u(r,x_3)}{r}\left(
                         \begin{array}{ccc}
                           -x_2 \\
                            x_1 \\
                            0
                         \end{array}
                       \right).
\end{equation}
Taking the fact $\text{div } U=0$ into account, it  follows   from  (\ref{d1.1}) that $u(r,x)$ satisfies
\begin{equation}\label{eq1.4}
-\Delta u+\frac{u}{r^2}+V(r,x_3)u=|u|^{p-2}u+|u|^4u,
\end{equation}
where $\Delta=\partial_r^2+\frac{1}{r}\partial_r+\partial_{x_3}^2$.
Correspondingly, we also have
\begin{equation}\label{eq1.00}
I(U)=J(u),
\end{equation}
where $J\in C^2(E,\R)$ is the energy functional of (\ref{eq1.4}) and given by
\begin{equation}\label{eq1.5}
J(u)=\frac{1}{2}\inte|\nabla u|^2+\frac{u^2}{r^2}+V(r,x_3)u^2dx-\frac{1}{p}\inte |u|^pdx-\frac{1}{6}\inte |u|^6dx, \ u\in E.
\end{equation}
Here, the space $E$ is defined by
\begin{equation*}
E:=\Big\{u(r,x_3)\in H^1(\R^3) \text{ and } \inte \frac{u^2}{r^2}dx<\infty\Big\},
\end{equation*}
and endowed with the norm of
$$\|u\|_E=\Big( \inte|\nabla u|^2+\frac{|u|^2}{r^2}+u^2dx\Big)^\frac{1}{2}.$$
Clearly $E\hookrightarrow H^1_s(\R^3):=\{u=u(r,x_3)\in H^1(\R^3)\}$. In general,  the integral  $\inte\frac{|u|^2}{r^2}dx$ cannot be controlled by $\inte|\nabla u|^2dx$ for all $u\in H^1_s(\R^3)$, so we have  $E\subsetneqq H_s^1(\R^3)$.

The derivative of $J(u)$ is given by
 \begin{equation*}
 \begin{split}
 \langle J'(u),\phi \rangle=\inte \nabla u\nabla \phi+\frac{u\phi}{r^2}+V(r,x_3)u\phi dx-\inte |u|^{p-2}u\phi dx-\inte u^5\phi dx, \ \forall\ \phi\in E.
 \end{split}
\end{equation*}
 We call $u\in E\setminus\{0\}$  a nontrivial weak solution of (\ref{eq1.4}) if $ \langle  J'(u),\phi\rangle=0$ for any $\phi\in E$, and denote the nontrivial solution set of (\ref{eq1.4}) by $$S=\big\{u: \text{ u is a nontrivial solution of (\ref{eq1.4})}\big\}.$$   Correspondingly, $u\in E$ is called  a  {\em ground state solution} of (\ref{eq1.4}) if
 $$u\in\big\{u \in S \text{ and } J(u)\leq J(v) \text{ for any } v\in S \big\}.$$

\begin{rem}
In view of the change  (\ref{eq1.04}), one can search for a special class of cylindrically symmetric solutions of (\ref{d1.1}) by studying equation (\ref{eq1.4}). Moreover, the ground state solution of (\ref{d1.1}) can be defined  similarly as that of (\ref{eq1.4}). As a consequence of (\ref{eq1.00}), any ground state solution of (\ref{eq1.4}) corresponds to a  ground state solution of (\ref{d1.1}) among all vector fields of type (\ref{eq1.04}). For these reasons, instead of search for solutions of equation (\ref{d1.1}) in the form of (\ref{eq1.04}),   we shall focus ourself on equation (\ref{eq1.4}) in what follows.
\end{rem}

Define
$$L^2_{\rm cyl}(\R^3):=\big\{u(r,x_3)\in L^2(\R^3) \big\} \text{
and }
W^{2,2}_{\rm cyl}(\R^3):=\big\{u(r,x_3)\in H^2(\R^3) \text{ and } \inte \frac{u^2}{r^2}dx<\infty\big\}.$$
Let $L: W^{2,2}_{\rm cyl}(\R^3)\subset L^2_{\rm cyl}(\R^3)\rightarrow L^2_{\rm cyl}(\R^3)$ be  given by $L:=-\Delta+\frac{1}{r^2}+V(r,x_3)$, which is  self-adjoint due to the fact that $\frac{1}{r^2}+V(r,x_3)$ is real valued, see e.g., \cite[Exampler p. 68]{T}. Assume that $V(r,x_3)$ satisfies
\begin{itemize}
  \item[(V)] $V(r,x_3)\in L^\infty(\R^3)$ is 1-periodic in $x_3$ direction  and  $0\not\in \sigma(L)$ where $\sigma$ denotes the spectrum of $L$ in $L^2_{\rm cyl}(\R^3)$.
\end{itemize}
Then our main result can be  stated as follows.
\begin{thm}\label{thm1}
Assume $V(r,x_3)$ satisfies $(V)$ and  $p\in(4,6)$, then (\ref{eq1.4}) has a ground state solution in $E$. Especially, if  $ \sigma(-\Delta+\frac{1}{r^2}+V)>0$ then (\ref{eq1.4}) possesses a ground state solution in $E$ for any $p\in(2,6)$.
\end{thm}

We intend to prove the above theorem by using Nehari-Pankov method as in \cite{BD}, where the subcritical case was studied. Let $\big(E(\lambda)\big)_{\lambda\in\R}$ be the spectral family of $L$. Set $E^-:= E(0)L^2_{\rm cyl}\cap E$ and $E^+:=(I-E(0))L^2_{\rm cyl}\cap E$. Then there exists an inner product $(\,,\, )$ such that the  corresponding norm $\| \cdot\|$ is equivalent  to $\|\cdot \|_E$, and for any $ u=u^++u^-\in E^++E^-$ there holds that 
$$\inte |\nabla u|^2+\frac{u^2}{r^2}+V(r,x_3)u^2dx= \|u^+\|^2-\|u^-\|^2.$$
We  note that if $\sigma (L)\subset(0,\infty)$ then $\text{dim }E^-=0$. Upon  the above decomposition for $E$,  functional $J$ defined in (\ref{eq1.5}) can be written as
\begin{equation*}
J(u)=\frac{1}{2}\|u^+\|^2-\frac{1}{2}\|u^-\|^2-\frac{1}{p}\inte |u|^pdx-\frac{1}{6}\inte |u|^6dx, \ u\in E,
\end{equation*}
and its derivative $J'(u)\in E^*$ is given by
\begin{equation}\label{eq1.12}
 \langle J'(u),\phi \rangle=( u^+,\phi^+)-( u^-,\phi^-)-\inte |u|^{p-2}u\phi dx-\inte u^5\phi dx, \ \forall\ \phi\in E.
\end{equation}
In section \ref{se2} we will consider the following minimization problem
\begin{equation} \label{eq1.11}
c:=\inf_{u\in \mathcal{N}} J(u),
\end{equation}
where  $\mathcal{N}$ is the  Nehari-Pankov manifold and  defined as
$$\mathcal{N}:=\big\{u\in E\setminus\{0\}: \langle J'(u), \phi\rangle=0 \ \forall\ \phi \in \R u\oplus E^- \big\}.$$
From (\ref{eq1.12})  one  see that $\langle J'(v),v \rangle<0$ if $v\in E^-\setminus\{0\}$, this indicates that $$\mathcal{N}\cap E^-=\emptyset.$$
If  we can prove that problem (\ref{eq1.11}) is attained by some $u\in E$, then it is a nontrivial solution, indeed, is a ground state solution of (\ref{eq1.4}) since all nontrivial solutions of (\ref{eq1.4})  belong to $\mathcal{N}$.  Theorem \ref{thm1} thus can be obtained.

We also note that  in view of   the existence of Sobolev critical exponent in equation (\ref{eq1.4}), to study problem (\ref{eq1.11}), in general, we  need to prove the energy of $J$ is strictly  less than some   critical value,  which is related to the following minimization problem
\begin{equation}\label{eq1.7}
\hat S:=\inf_{0\not=u\in E}\frac{\inte |\nabla u|^2+\frac{|u|^2}{r^2}dx}{(\inte |u|^6 dx)^\frac{1}{3}}.
\end{equation}
Repeating the argument of \cite[Theorem 1.3]{CSW}, one can prove that (\ref{eq1.7}) can be achieved by a nonnegative minimizer $\Phi(r,x_3)\in E$, which is  Steiner-symmetric with respect to $x_3$. Moreover, $\Phi$ satisfies the following equation
\begin{equation}\label{eq1.8}
-\Delta \Phi+\frac{\Phi}{r^2}=\Phi^5 \ \text{ in }\ \R^3.
\end{equation}
From (\ref{eq1.7}) and (\ref{eq1.8}) one can easily check that
\begin{equation}\label{eq1.9}
\inte |\nabla \Phi|^2+\frac{|\Phi|^2}{r^2}dx=\inte|\Phi|^6dx=\hat S^\frac{3}{2}.
\end{equation}
Also,  the results of \cite[propositons 5 and 6]{BBR} indicate that  $\Phi\in L^\infty(\R^3)$ and
\begin{equation}\label{eq1.10}
\limsup_{|x|\to\infty}|x|^\nu\Phi(x)<\infty \  \text{ for any } \upsilon<\frac{1+\sqrt{5}}{2}.
\end{equation}
This especially yields  that $\Phi\in L^2(\R^3)\cap L^6(\R^3)$.

In what follows we denote by $C$ a positive constant, which may be different from line to line. The open ball centered at $x_0$ with radius $R$ is denoted by $B_R(x_0)$ and  $|\cdot|_{L^q}$ means the norm in $L^q(\R^3)$ .

\section{Proof of Theorem \ref{thm1}}\label{se2}
In this section we give the proof of Theorem \ref{thm1} by using the Nehari-Pankov method.  For this purpose, we recall two key lemmas which have been  proved in previous papers, so we omit the proof here. Firstly, similar to the arguments of Propositions 2.2 and 2.4 in \cite{CS}, we  have the following lemma.

\begin{lem}\label{le2.1}
If $V(x)$ satisfies $(V)$, then for any $\mu\in\R$ there exist positive constants  $C_0$ and $C_1$ such that $|u|_{L^\infty}\le C_0\|u\|_E \le C_1 |u|_{L^2}$ for all $u\in E(\mu)L^2_{\rm cyl}(\R^3)\cap E$.
\end{lem}

 Repeating the  arguments of \cite[lemmas 19-21]{BD} one can prove  that for any $u\in \mathcal{N}$ the norm of $\nabla J(u)$ can be controlled by its tangential component. Here $\nabla J(u)$ denotes the unique correspondence of $J'(u)$ in $E$. Precisely, we have
\begin{lem}\label{le2.4}
Let $\mathcal{N}_0 $ be a bounded subset of $\mathcal{N}$. There exists $C_0>0$ such that the following holds: for any   $u\in \mathcal{N}_0$, $\nabla J(u)=\tau+\sigma$ where $\tau\in T_u\mathcal{N}$ is the tangential component of $\nabla J(u)$ and $\sigma\bot \tau$ is the transversal component of $\nabla J(u)$, then,
$$\|\nabla J(u)\|\leq C_0\|\tau\|.$$
\end{lem}

\vskip 0.1truein
The following lemma tells that  the norm $\|u\|$ and $J(u)$ are both bounded below  on the manifold $\mathcal{N}$.

\begin{lem}\label{le2.30}
Assume  $V(x)$ satisfies $(V)$, then there exist $C_1,C_2>0$ such that
\begin{equation}\label{eq2.100}
\|u\|\geq C_1 \ \text{ and } \|u\|\leq C_2(J(u))^\frac{5}{6} \ \text{ for all } u\in\mathcal{N}.
\end{equation}
\end{lem}
\noindent\textbf{Proof.} For any $u\in \mathcal{N}$, then
\begin{equation}\label{eq2.0}J(u)=\big(\frac{1}{2}-\frac{1}{p}\big)\inte|u|^pdx+\frac{1}{3}\inte|u|^6dx\geq0.\end{equation}
Therefore, we have $\|u^-\|\leq\|u^+\|$.
Moreover, it follows from  $\langle J'(u), u^+\rangle=0$ that for any $\delta>0$,
\begin{align}
\frac{\|u\|^2}{2}&\le\|u^+\|^2=\inte|u|^{p-2}uu^+dx+\inte u^5u^+dx\leq \delta\inte |u||u^+|dx+C|u|^5|u^+|dx\nonumber\\
&\leq \delta|u|_{L^2}|u^+|_{L^2}+C(\delta)|u|_{L^6}^5|u^+|_{L^6}\leq \delta\|u\|^2+C\|u\|^6. \label{eq2.01}
\end{align}
Taking $\delta<\frac{1}{4}$, we then deduce that  there exists $C_1>0$ such that
$
\|u\|\ge C_1>0
$.
Furthermore, it follows from (\ref{eq2.0}) and (\ref{eq2.01}) that
\begin{align*}
\frac{\|u\|^2}{2}&\leq \delta\|u\|^2+C|u|_{L^6}^5|u^+|_{L^6}\leq \delta\|u\|^2+C|u|_{L^6}^5\|u\| \\
&\leq \delta\|u\|^2+C(J(u))^\frac{5}{6}\|u\|.
\end{align*}
This indicates that $(J(u))^\frac{5}{6}\geq\|u\|\geq C_1>0 $.
\qed

Stimulated  by  the arguments of \cite{BBR,BGR,CS,SZ}, we can  apply   Lemmas \ref{le2.1} - \ref{le2.4} to obtain the following theorem, which yields Theorem \ref{thm1} if  condition  $c<\frac{1}{3}\hat S^\frac{3}{2}$ is  satisfied.
\begin{thm}\label{thm2}
If $c<\frac{1}{3}\hat S^\frac{3}{2}$ and $p\in(2,6)$, then the minimum in (\ref{eq1.11}) can be achieved by some $u\in E$. Moreover, $u$ is a ground state solution of (\ref{eq1.4}).
\end{thm}

\noindent\textbf{Proof. } Using Ekeland's variational principle \cite{St}, there exists a minimizing sequence $\{u_n\}\subset \mathcal{N}$ of $c$  such that $J(u_n)\overset n\to c$ and $(J|_{\mathcal{N}})'(u_n)\overset{n}\to 0$. Furthermore, from  Lemma \ref{le2.30} we see that $\{u_n\}$ is  bounded in $E$. Thus,  $J'(u_n)\overset n\to0$ by applying Lemma \ref{le2.4}. In summary, we obtained  a  Palais-Smale sequence
 \begin{equation}\label{eq2.9}
J(u_n)\overset{n}\to c\ \text{ and } J'(u_n)\overset{n}\to0 \text{ in } E^*.
\end{equation}

We first rule out the case of {\em vanishing}. Assume for contradiction that
   \begin{equation}\label{eq2.15}
   \lim_{n\to\infty}\sup_{y\in\R^3}\int_{B_R(y)}|u_n|^2dx=0 \text{ for any }R>0.
   \end{equation}
It  follows from  \cite[Lemma I.1]{l2} that
\begin{equation}\label{eq2.11}
u_n\overset{n}\to0 \text{ in }L^q(\R^3)\text{ for any }q\in(2,6).\end{equation}
Thus,
\begin{equation}
J(u_n)-\frac{1}{2}\langle J'(u_n), u_n\rangle=\frac{1}{3}\inte|u_n|^6dx+o(1)\overset{n}\to c.
\end{equation}

From  (\ref{eq2.11}) and  Lemma \ref{le2.1} we deduce from the fact $\langle J'(u_n), u_n^-\rangle=0$ that
\begin{equation}
\begin{split}\label{eq2.12}
\|u_n^-\|^2&=-\inte |u_n|^4u_nu_n^-dx-\inte|u_n|^{p-2}u_nu_n^-dx\\
&\leq |u_n|_{L^5}^5|u_n^-|_{L^\infty}+|u_n|_{L^p}^{p-1}|u_n^-|_{L^p}+o(1)\overset{n}\to0.
\end{split}
\end{equation}
Rewrite  $u_n^+=w_n+z_n$, where $w_n\in E(\mu)L^2_{\rm cyl}\cap E$ and $z_n\in (I-E(\mu))L^2_{\rm cyl}\cap E$, $\mu>0$ is larger enough (which will be determined later).  As a consequence of (\ref{eq2.9}) we have  $\langle J'(u_n), w_n\rangle=o(1)$.  Similar to the argument of (\ref{eq2.12}), one can obtain that
 $$\|w_n\|^2=\inte |u_n|^4u_nw_ndx+\inte|u_n|^{p-2}u_nw_ndx+o(1)\overset{n}\to0.$$
 Therefore,
 \begin{equation}
 |u_n-z_n|_{L^6}=|u_n^-+w_n|_{L^6}\leq C\|u_n^-+w_n\|\overset{n}\to0,
 \end{equation}
 and
 \begin{equation}
 \begin{split}
\|z_n\|^2&=\inte |\nabla z_n|^2+\frac{|z_n|^2}{r^2}+V(r,x_3)|z_n|^2dx\\
&=\langle J'(u_n), u_n\rangle-\|w_n\|^2+\|u_n^-\|^2+\inte |u_n|^6dx+\inte|u_n|^pdx\\
&=\inte |u_n|^6dx+o(1).
\end{split}
\end{equation}
Since $\mu\inte|z_n|^2dx\leq \inte |\nabla z_n|^2+\frac{|z_n|^2}{r^2}+V|z_n|^2dx$, we have $(\mu-|V|_{L^\infty})\inte|z_n|^2dx\leq\inte |\nabla z_n|^2+\frac{|z_n|^2}{r^2}dx $. For any $\delta>0$, we take $\mu>0$ large enough such that $|V|_{L^\infty}\leq \delta(\mu-|V|_{L^\infty})$. Then
$$\delta\inte |\nabla z_n|^2+\frac{|z_n|^2}{r^2}dx\geq \delta(\mu-|V|_{L^\infty})\inte|z_n|^2dx\geq |V|_{L^\infty} \inte|z_n|^2dx.$$
This implies that
\begin{equation}\label{eq2.211}(1-\delta)\inte|\nabla z_n|^2+\frac{|z_n|^2}{r^2}dx\leq\inte |\nabla z_n|^2+\frac{|z_n|^2}{r^2}+V|z_n|^2dx=\|z_n\|^2.\end{equation}
We deduce from (\ref{eq2.11})-(\ref{eq2.211}) that
\begin{equation*}
\begin{split}
c+o(1)&=J(u_n)=\frac{1}{2}\|z_n\|^2+\frac{1}{2}\|w_n\|^2-\frac{1}{2}\|u_n^-\|^2-\frac{1}{6}\inte|u_n|^6dx+\frac{1}{p}\inte|u_n|^pdx\\
&=\frac{1}{2}\|z_n\|^2-\frac{1}{6}\inte|u_n|^6dx+o(1)\geq\frac{1-\delta}{2}\inte|\nabla z_n|^2+\frac{|z_n|^2}{r^2}dx-\frac{c}{2}+o(1)\\
&\geq\frac{(1-\delta)\hat S}{2}|z_n|_{L^6}^2-\frac{c}{2}+o(1)=\frac{(1-\delta)\hat S}{2}|u_n|_{L^6}^2-\frac{c}{2}+o(1)\\
&=\frac{(1-\delta)(3c)^\frac{1}{3}\hat S}{2}-\frac{c}{2}+o(1).
\end{split}
\end{equation*}
From (\ref{eq2.100}) we see that $c>0$.  Letting $n\to\infty$, the above inequality thus implies that
$c\geq\frac{1}{3}(1-\delta)^\frac{3}{2}\hat S^\frac{3}{2}$ for all $\delta>0$. This leads to a contradiction  for the condition $c<\frac{1}{3}\hat S^\frac{3}{2}$  is assumed.

Now we have proved that {\em vanishing} (\ref{eq2.15}) cannot occur. Thus, there exist $R,\eta>0$ and a sequence $\{x_n=(y_n,x_{3n})\}\subset\R^3$ (without loss of generality, we assume that $x_n\in \mathbb{Z}^3$) such that
\begin{equation}\label{eq2.16}
\limsup_{n\to\infty}\int_{B_R(x_n)}|u_n|^2dx\geq\eta>0.
\end{equation}
We claim that
\begin{equation}\label{eq2.17}
\limsup_{n\to\infty}|y_n|<\infty.\end{equation}
Firstly,  since  $u_n(x)=u_n(|y|,x_3)$ is cylindrically symmetric, thus
\begin{equation}\label{eq2.18}
\int_{B_R(x_n)}|u_n|^2dx=\int_{B_R(g\cdot y_n,x_{3n})}|u_n|^2dx \ \text{ for all }\ g\in O(2),
\end{equation}
where  $O(2)$ is the  orthogonal group of $\R^2$.
On the contrary, if (\ref{eq2.17}) is false, i.e.,  there exists a subsequence, still denoted by $n$, such that $\lim_{n\to\infty}|y_n|=\infty$.  For the point $\bar x=(\bar y,x_{3n})$ satisfying  $|\bar y|=|y_n|$, it follows  from
\cite[proposition 17]{BBR} or the proving of \cite[Theorem 1.3]{BD} that, the number of disjoint balls $B_R(\bar x)$ goes to infinity if $\lim_{n\to\infty}|y_n|=\infty$. This is impossible since $\{|u_n|_{L^2}\}$ is bounded.  Claim (\ref{eq2.17}) is thus  proved.


Let $\bar x_n=(0, x_{3n})$, in view of (\ref{eq2.17}) we may assume that (\ref{eq2.16}) still holds by replacing $x_n$ with $\bar x_n$. Set $\bar u_n(x)=u_n(x+\bar x_n)$, then
\begin{equation}\label{eq2.19}
\limsup_{n\to\infty}\int_{B_R(0)}|\bar u_n|^2dx\geq\eta>0.
\end{equation}
Moreover, since $V(r,x_3)$ is 1-periodic in $x_3$ direction, we see that  $\{\bar u_n\}\subset \mathcal{N}$ satisfies
 \begin{equation}\label{eq2.90}
J(u_n)=J(\bar u_n)\overset{n}\to c\ \text{ and } J'(u_n)=J'(\bar u_n)\overset{n}\to0 \text{ in } E^*.
\end{equation}
 Hence,  by passing to subsequence, there exists $u\in E$ such that
$\bar u_n\overset{n}\rightharpoonup u$ in $E$ and $J'(u)=0$. From (\ref{eq2.19}) we see that $u\not\equiv0$, so $u$ is a nontrivial solution of (\ref{eq1.4}) by (\ref{eq2.90}). Furthermore, in view of $u,\bar u_n \in \mathcal{N}$, it follows from Fatou's lemma  that
\begin{align*}
c&=\lim_{n\to\infty}J(\bar u_n)=\lim_{n\to\infty}\Big[\big(\frac{1}{2}-\frac{1}{p}\big)\inte|\bar u_n|^pdx+\frac{1}{3}\inte|\bar u_n|^6dx\Big]\\
&\geq\big(\frac{1}{2}-\frac{1}{p}\big)\inte| u|^pdx+\frac{1}{3}\inte| u|^6dx= J(u)\geq c.
\end{align*}
This implies that $u$ is a minimizer of (\ref{eq1.11}), so it is a ground state solution of (\ref{eq1.4}).
\qed
\vskip 0.2 truein
From Theorem \ref{thm2} we see that, to finish the proof of Theorem \ref{thm1} it remains to show that the energy value $c<\frac{1}{3}\hat S^\frac{3}{2}$. Our following three lemmas tell that this condition can be verified under the assumptions of Theorem \ref{thm1}.

\begin{lem}\label{le2.2}
Let $\Phi$ be the nonnegative solution given by equation (\ref{eq1.8})  and set
\begin{equation*}
\varphi_\varepsilon(x)=\varepsilon^{-\frac{1}{2}}\Phi(\frac{x}{\eps})=\varphi_\varepsilon^+(x)+\varphi_\varepsilon^-(x)\in E^+\oplus E^-,
\end{equation*}
then,\begin{equation}\label{eq2.3}
     \inte|\nabla \varphi_\varepsilon|^2+\frac{|\varphi_\varepsilon|^2}{r^2}dx=\inte |\varphi_\varepsilon|^6=\hat S^\frac{3}{2},  \inte |\varphi_\varepsilon|^qdx=O(\eps^{-\frac{q}{2}+3}) \ \forall\ q\in[2,6);\end{equation}
\begin{equation}\label{eq2.4}
       |\varphi_\eps^-|_{L^\infty},\|\varphi_\eps^-\|, |\varphi_\eps^-|_{L^2}\leq O(\eps);
       \end{equation}
\begin{equation} \label{eq2.44}
|\varphi_\varepsilon^+|_{L^5}^5\leq C \eps^\frac{1}{2},\ \   \|\varphi_\varepsilon^+\|^2= \hat S^\frac{3}{2}+ O(\eps^2);\end{equation}
and
 \begin{equation}\label{eq2.60}      \Big|\inte|\varphi_\varepsilon|^6dx-\inte|\varphi_\varepsilon^+|^6dx\Big|\leq C \eps^\frac{3}{2}.\end{equation}
\end{lem}

\noindent\textbf{Proof.} From (\ref{eq1.9}) and (\ref{eq1.10}) we see that  $\Phi(x)\in L^2(\R^3)\cap L^6(\R^3)$, then direct calculations show that  (\ref{eq2.3}) holds. Using (\ref{eq2.3}) we see that $|\varphi_\eps^-|_{L^2}\leq C|\varphi_\eps|_{L^2}\leq C\eps$, this yields (\ref{eq2.4}) by applying Lemma \ref{le2.1}.
 From (\ref{eq2.3}) and (\ref{eq2.4}), we have
 $$\inte |\varphi_\varepsilon|^5dx\leq C\eps^{\frac{1}{2}}\text{ and } \inte |\varphi_\varepsilon^-|^5dx\leq |\varphi_\varepsilon^-|_{L^\infty}^3|\varphi_\varepsilon^-|_{L^2}^2\leq C\eps^{5}.$$
 This indicates that
$$|\varphi_\eps^+|_{L^5}^5=|\varphi_\eps-\varphi_\eps^-|_{L^5}^5\leq C(|\varphi_\eps|_{L^5}^5+|\varphi_\eps^-|_{L^5}^5)\leq C\eps^\frac{1}{2}.$$
Morever, in view of $V\in L^\infty(\R^3)$ there holds that
$$\|\varphi_\varepsilon^+\|^2=\inte|\nabla \varphi_\varepsilon|^2+\frac{\varphi_\varepsilon^2}{r^2}+V|\varphi_\eps|^2dx+\|\varphi_\eps^-\|^2=\hat S^\frac{3}{2}+O(\eps^2).$$
Thus (\ref{eq2.44}) is obtained.

 To prove (\ref{eq2.60}), we first deduce from the convexity of $|\cdot|_{L^6}$ and Lemma \ref{le2.1} that
\begin{equation}\label{eq2.61}
\begin{split}
\inte|\varphi_\eps|^6dx&=\inte|\varphi_\eps^++\varphi_\eps^-|^6dx\geq \inte|\varphi_\eps^+|^6dx+6\inte(\varphi_\eps^+)^5\varphi_\eps^-dx\\
&\geq\inte|\varphi_\eps^+|^6dx-6\inte|\varphi_\eps^+|^5dx|\varphi_\eps^-|_{L^\infty}\geq\inte|\varphi_\eps^+|^6dx-C\eps^\frac{3}{2},
\end{split}
\end{equation}
where (\ref{eq2.4}) and (\ref{eq2.44}) is used in the last inequality. Similarly, it follows from (\ref{eq2.3}) and (\ref{eq2.4}) that
$$\inte|\varphi_\eps^+|^6dx=\inte|\varphi_\eps-\varphi_\eps^-|^6dx\geq \inte|\varphi_\eps|^6dx-C\eps^\frac{3}{2}.$$
This together with (\ref{eq2.61}) yields (\ref{eq2.60}). \qed

\vskip 0.2truein

Let
\begin{equation*}
M_\eps:=\big\{u=u^-+t\varphi_\eps^+: u^-\in E^- \text{ and } \|u\|<R\big\},
\end{equation*}
Then we have

\begin{lem}\label{le2.3}
\begin{itemize}
  \item [\rm(i)] If $\varepsilon>0$ is small enough,  there exists $R>0$ independent of $\eps$ such that
  \begin{equation}\label{eq2.5}
  \sup_{u\in \partial M_\eps}J(u)\leq 0 \text{ and } \sup_{u\in M_\eps}J(u)<\infty.\end{equation}
      \item [\rm(ii)] There exist $0<\rho<R$ and $\alpha>0$ such that $J(u)\geq\alpha$ for all $u\in E^+\cap \partial B_{\rho}(0)$.
\end{itemize}

\end{lem}

\noindent\textbf{Proof. } \textbf{(i).} To prove (\ref{eq2.5}) it is sufficient to prove that there exists $R>0$ large enough such that
\begin{equation}\label{eq2.6}
J(u^-+t\varphi_\eps^+)\leq0  \ \ \text{ for any  }\|u^-+t\varphi_\eps^+\|\geq R.
\end{equation}
We argue  in two different cases. Firstly, if $\|t\varphi_\eps^+\|\leq \frac{1}{2}\|u^-\|$,
\begin{equation}\label{eq2.7}
J(u^-+t\varphi_\eps^+)\leq \frac{1}{2}\|t\varphi_\eps^+\|^2-\frac{1}{2}\|u^-\|^2\leq -\frac{1}{4}\|u^-\|^2<0.
\end{equation}
On the other hand, if $\|t\varphi_\eps^+\|>\frac{1}{2}\|u^-\|$, then
$$R^2\leq \|u^-+t\varphi_\eps^+\|^2=\|u^-\|^2+\|t\varphi_\eps^+\|^2\leq 5t^2\|\varphi_\eps^+\|^2,$$
together with (\ref{eq2.44}) we see that if $\eps>0$ is small,
\begin{equation}\label{eq2.10}t^2\geq\frac{R^2}{5\|\varphi_\eps^+\|^2}\geq\frac{R^2}{5(\hat S^\frac{3}{2}+1)}.\end{equation}
 Moreover, it follows from Lemma \ref{le2.1}, (\ref{eq2.3}), (\ref{eq2.44}) and (\ref{eq2.60}) that if $\eps>0$ is small enough,
\begin{align}
&\inte |u^-+t\varphi_\eps^+|^6dx\geq t^6\inte |\varphi_\eps^+|^6dx+6\inte (t\varphi_\eps^+)^5u^-dx \nonumber\\
&\geq t^6\inte |\varphi_\eps^+|^6dx-6|t|^5\|u^-\|\inte |\varphi_\eps^+|^5dx\geq t^6\inte |\varphi_\eps^+|^6dx-C\eps^\frac{1}{2}|t|^5\|u^-\|\label{eq2.8} \\
&\geq  t^6\inte |\varphi_\eps^+|^6dx-C\eps^\frac{1}{2}t^6\|\varphi_\eps^+\|\geq t^6(\hat S^\frac{3}{2}+O(\eps^\frac{3}{2}))-C\hat S^\frac{3}{4}\eps^\frac{1}{2}t^6\geq \frac{\hat S^\frac{3}{2}}{2}t^6.\nonumber
\end{align}
 Combine with (\ref{eq2.10}), we see that if $R>0$ is suitable large,
\begin{equation*}
\begin{split}J(u^-+t\varphi_\eps^+)&\leq\frac{1}{2}\|t\varphi_\eps^+\|^2-\frac{1}{6}\inte |u^-+t\varphi_\eps^+|^6dx\leq \hat S^\frac{3}{2}t^2-\frac{\hat S^\frac{3}{2}}{12}t^6<0.
\end{split}\end{equation*}
This together with (\ref{eq2.7}) gives (\ref{eq2.6}). From (\ref{eq2.6}) we see that $\sup_{u\in \partial M_\eps}J(u)\leq 0 $. Moreover, since $J$ maps bounded sets into bounded sets, therefore  $\sup_{u\in M_\eps}J(u)<\infty$.

\textbf{(ii).} Note that  $|u|_{L^q}\leq C\|u\| $ for any $q\in[2,6]$. For any  $u\in E^+\cap \partial B_{\rho}(0)$ with $\rho$ small enough,  we have
$$J(u)\geq \frac{1}{2}\|u\|^2-C_1\|u\|^p-C_2\|u\|^6\geq C_3\|u\|^2-C_4\|u\|^6>C_\rho>0. $$
 \qed

\begin{lem}\label{le2.5}
Let $c$ be defined as in (\ref{eq1.11}), then \begin{equation}\label{eq2.200}
c<\frac{1}{3}\hat S^\frac{3}{2},\end{equation}
if one of the following conditions holds:
\begin{itemize}
  \item [\rm(i)] $V(r,x_3)$ satisfies $(V)$ and $p\in(4,6)$;
  \item [\rm(ii)] $V(r,x_3)\in L^\infty(\R^3)$ satisfies $0< \sigma(-\Delta+\frac{1}{r^2}+V)$ and $p\in(2,6)$.
\end{itemize}

\end{lem}
\noindent\textbf{Proof. }
Set$$\gamma:=\sup_{v\in \R \varphi_\eps^+\oplus E^-}J(v),$$
From the proof of Lemma \ref{le2.3} (i), we see that $J(v)\leq0$ if $v\in M_\eps^c$. This  together with
 Lemma \ref{le2.3} (ii) indicates that
 $$\gamma=\sup_{v\in M_\eps}J(v)\geq \alpha>0.$$
Therefore, there exists $\{u_n=t_n\varphi_\eps^++u_n^-\}\subset M_\eps$ such that $\lim_{n\to\infty}J(u_n)=\gamma$.   Up to subsequence, we can assume that there exists $u_0= t_0\varphi_\eps^++u_0^-\in \overline M_\eps$ such that $$u_n=t_n\varphi_\eps^++u_n^-\overset{n}\rightharpoonup u_0=t_0\varphi_\eps^++u_0^- \ \text{ in } E.$$
This implies that $t_n\overset{n}\to t_0$. Moreover, since $\psi(u):=\frac{1}{2}\|u^-\|^2+\frac{1}{6}\inte |u|^6dx+\frac{1}{p}\inte |u|^pdx$ is weakly sequentially lower semicontinuous in $E$, we thus deduce that
$$\gamma=\lim_{n\to\infty}J(u_n)\le J(u_0)\leq \gamma.$$
Thus, $u_0$ is a maximum point of $J$ on $\R \varphi_\eps^+\oplus E^-$. This yields that $u_0\in \mathcal{N}$ and
$c\leq \sup_{v\in M_\eps}J(v).$ Hence, to obtain (\ref{eq2.200}), it is sufficient to prove that
 \begin{equation}\label{eq2.277}
\sup_{v\in M_\eps}J(v)<\frac{1}{3}\hat S^\frac{3}{2}.
\end{equation}

\textbf{(i).}
For any $u=u^-+t\varphi_\eps^+\in M_\eps$, since $\|u\|\leq R$, we deduce from  (\ref{eq2.44}) that if $\eps>0$ is small enough, then there exists $\beta=\beta(R)>0$ such that
\begin{equation}\label{eq2.28}|t|\le\beta \text{ and }\|u^-\|\leq \beta.\end{equation}
 It then follows from (\ref{eq2.44}) that
\begin{equation}
\begin{split}\label{eq2.21}
&\inte |\nabla u|^2+\frac{|u|^2}{r^2}+V|u|^2dx=\|t\varphi_\eps^+\|^2-\|u^-\|^2=\hat S^\frac{3}{2}t^2+O(\eps^2)-\|u^-\|^2.
\end{split}
\end{equation}
 Similar to (\ref{eq2.8}), we deduce from Lemmas \ref{le2.1}, \ref{le2.2} and (\ref{eq2.28}) that
\begin{align}
&\inte |u|^6dx\geq t^6\inte |\varphi_\eps^+|^6dx-C\eps^\frac{1}{2}|t|^5\|u^-\|\geq t^6\inte |\varphi_\eps|^6dx-C\eps^\frac{3}{2}-C\beta^5\eps^\frac{1}{2}\|u^-\|\nonumber\\
&=\hat S^\frac{3}{2}t^6-C_1\eps^\frac{3}{2}-C_1\eps^\frac{1}{2}\|u^-\|,\label{eq2.22}
\end{align}
and
\begin{align}
&\inte |u|^pdx=\inte|t\varphi_\eps-t\varphi_\eps^-+u^-|^pdx\geq |t|^p\inte|\varphi_\eps|^pdx-p\inte|t\varphi_\eps|^{p-1}|t\varphi_\eps^-+u^-|dx\nonumber\\
&\geq |t|^p\inte|\varphi_\eps|^pdx-C(\beta|\varphi_\eps^-|_{L^\infty}+|u^-|_{L^\infty})|t|^{p-1}\inte|\varphi_\eps|^{p-1}dx\nonumber\\
&\geq  |t|^p\inte|\varphi_\eps|^pdx-C(\beta|\varphi_\eps^-|_{L^\infty}+\|u^-\|)\beta^{p-1}\inte|\varphi_\eps|^{p-1}dx\nonumber\\
&\geq |t|^p\inte|\varphi_\eps|^pdx-C_2\inte|t\varphi_\eps|^{p-1}dx\geq |t|^pC_2\eps^{-\frac{p}{2}+3}-C_2\eps^{-\frac{p-1}{2}+3}.\label{eq2.23}
\end{align}
Using (\ref{eq2.21})-(\ref{eq2.23}) and the fact that  $-\|u^-\|^2+C_1\eps^\frac{1}{2}\|u^-\|\leq C_3\eps$, we see that
\begin{equation}\label{eq2.24}
J(u)\leq \hat S^\frac{3}{2}\big(\frac{t^2}{2}-\frac{t^6}{6}\big)-|t|^p C_2\eps^{-\frac{p}{2}+3}+C_2\eps^{-\frac{p-1}{2}+3}+C_3\eps.
\end{equation}
If $ |t|\leq \frac{1}{2}$, then $\big(\frac{t^2}{2}-\frac{t^6}{6}\big)< \frac{1}{8}$, and it is obvious to see that if $\eps>0$ is small enough,
 \begin{equation}\label{eq2.26}J(u)\leq \frac{1}{8}\hat S^\frac{3}{2}+C\eps^{-\frac{p}{2}+3}+C_3\eps< \frac{1}{4}\hat S^\frac{3}{2}.\end{equation}
 On the other hand, if $\frac12<|t|\leq \beta$, we then deduce from (\ref{eq2.24}) and the fact $\sup_{t\in\R}\big(\frac{t^2}{2}-\frac{t^6}{6}\big)=\frac{1}{3}$ that
 \begin{equation}\label{eq2.27}J(u)\leq \frac{1}{3}\hat S^\frac{3}{2}-\frac{1}{2^p}C_2\eps^{-\frac{p}{2}+3}+C_2\eps^{-\frac{p-1}{2}+3}+C_3\eps\leq \frac{1}{3}\hat S^\frac{3}{2}-\frac{1}{2^p}C_4\eps^{-\frac{p}{2}+3},\end{equation}
where $p\in(4,6)$ is used in the last inequality. In brief, we deduce from  (\ref{eq2.26}) and (\ref{eq2.27}) that
$$J(u)\leq \max\big\{\frac{1}{4}\hat S^\frac{3}{2},\frac{1}{3}\hat S^\frac{3}{2}-\frac{1}{2^p}C\eps^{-\frac{p}{2}+3}\big\} \text{ for all } u\in M_\eps.$$
Thus (\ref{eq2.277}) holds if $\eps>0$ is small.

 \textbf{(ii).} If $0< \sigma(-\Delta+\frac{1}{r^2}+V)$, then $E^-=\{0\}$ and thus
$M_\eps=\big\{t\varphi_\eps:  \|t\varphi_\eps\|<R\big\}$.
It then follows from Lemma \ref{le2.2} and $p\in(2,6)$ that
\begin{align*}
J(t\varphi_\eps)&=\frac{t^2}{2}\inte |\nabla \varphi_\eps|^2+\frac{|\varphi_\eps|^2}{r^2}+V|\varphi_\eps|^2dx-\frac{|t|^p}{p}\inte |\varphi_\eps|^pdx-\frac{t^6}{6}\inte |\varphi_\eps|^6dx\\
&=\hat S^\frac{3}{2}\big(\frac{t^2}{2}-\frac{t^6}{6}\big)+C_5\eps^2-\frac{|t|^p}{p}C_6\eps^{-\frac{p}{2}+3}\leq \frac{1}{3}\hat S^\frac{3}{2}-C_7\eps^{-\frac{p}{2}+3}.
\end{align*}
Then (\ref{eq2.277}) follows by taking $\eps>0$ small enough.\qed

\vskip 0.2 truein

\noindent {\bf Acknowledgements:} The author would like to thank the referees for the very helpful comments upon which the paper is revised. This research was supported by
    the National Natural Science Foundation of China under grant Nos. 11501555, 11471330 and 11471331.

\end{document}